\documentclass[12pt]{article}
\usepackage{lmodern}
\usepackage[margin=1in]{geometry}
\usepackage{amsmath,amssymb,amsthm}
\usepackage{enumitem}
\usepackage{changepage}
\usepackage{pgfplots}
\usepackage{url}
\usepackage{subcaption}
\usepackage{fancyhdr}
\usepackage{newtxtext}
\usepackage{comment}
\usepackage{lettrine}
\usepackage{float}
\usepackage{fancybox,framed}
\theoremstyle{definition}
\newtheorem{definition}{Definition}
\newtheorem{theorem}{Theorem}
\newtheorem{proposition}{Proposition}
\newtheorem{lemma}{Lemma}
\newtheorem{corollary}{Corollary}[proposition]
\renewcommand\thecorollary{\theproposition.\arabic{corollary}}
\newcommand{\Z}{\mathbb Z}

\pgfplotsset{compat=1.12}

\title{Iteration Sums of The Euler's Totient Function Regarding Powers of Fermat Primes}
\author{
Xiang Li, Allison Pacelli \\
Pioneer Research Number Theory
}
\date{October 9 2025}

\newcounter{prob}

\begin{document}
\maketitle
\pagestyle{fancy}
\fancyhead{} 
\fancyhead[L]{Pioneer Research}
\fancyhead[C]{\textbf{The Totient Function}}
\fancyhead[R]{Xiang Li}

\begin{abstract}
 Euler’s totient function $\phi(n)$ plays a central role in number theory and is applied in areas such as cryptography. In this paper, we study iterations of the totient function. We first prove that for any integer $n>2$, iteratively applying $\phi$ eventually yields the value $2$. Motivated by this terminal behavior, we examine sums of iterated totient values of the form $\phi(n)+\phi(\phi(n))+\phi(\phi(\phi(n)))+\cdots+\phi(2)$, where the summation terminates at $\phi(2)$. We show that for all integers of the form $n = 3^k$, this sum is equal to $n$. We then extend this result to all powers of Fermat primes, deriving a closed-form expression for the corresponding summations.
\end{abstract}

\section{Introduction}

\lettrine{P}{rimes, Divisibility, Factorization, and Modular Arithmetic} have long been the central topics of number theory. As one of the oldest branches of mathematics, studies on these concepts have begun from ancient times in Egypt and Greece \cite{Number_theory_oldest_branches_in_mathematics} up to modern times, with numerous scholars investing their time and genius. Many problems were solved, while more problems were proposed. Assisting with proofs, functions served as a major tool in mathematics. The Euler’s totient function, intertwined with multiple concepts of number theory, is the major concern of this paper. The word ``totient'' originates from the Latin word ``\textit{tot}", which means ``so many", as an answer to the question ``\textit{Quot}"(``How many") \cite{totient_in_Latin}. In other literature, the Toitient function is also referred to as the Phi function, because it is denoted by $\phi(n)$. The totient function counts the number of positive integers relatively prime to $n$ that are less than or equal to $n$.\\

When there exists a function, then naturally, the question of applying the function more than once comes to mind. When we iterate a function multiple times, each time the previous output of the function becomes the new input of the function, which again produces an output that could be used for the next iteration. Therefore, the totient function suites well to iterate, given that both its inputs and outputs are positive integers. Every output in the range would be a valid input in the domain of the totient function, ensuring that undefined values never arise. This paper provides a new approach to sum up iterations of the totient function: $\phi(n)+\phi(\phi(n))+\phi(\phi(\phi(n)))+\cdots$. With this approach, we prove a formula that gives the iterated totient sums of all powers of Fermat primes, which produced a clean corollary for all powers of 3.\\

These investigations are important because, first of all, the study of the totient function impacts various disciplines in mathematics. For example, multiplicative functions are introduced when the totient function is proved to be multiplicative, but not completely multiplicative. The Chinese remainder theorem plays a role in the proof of The totient function's multiplicativity, also affecting many studies of modular arithmetic, such as linear congruences. Then, the iteration sums of the totient function bring the topic of Fermat Primes, a special type of number with unsolved conjectures. In fact, the totient function itself has many well known conjectures such as Carmichael's Conjecture \cite{Open_Conjecture_Carmichael} and Lehmer's Problem \cite{Open_Conjecture_Lehmer}. Although this paper focuses on theory instead of application, there are many direct applications of number theory in real life. With regards to the context of this paper, the Chinese Remainder Theorem and the totient function are utilized in cryptography \cite{chinese_remainder_thm}. Specifically, Euler's totient Theorem provided the algorithm for one of the most used public key cryptosystems in the world--- RSA \cite{RSA}. By formalizing patterns into theorems, number theory has laid the foundation for transformative technology advancements.

\section{Historical Background}

In 1763, the first definition of the totient function appeared in Euler's work, \textit{Theoremata arithmetica nova methodo demonstrata} \cite{euler_totient}. However, Euler did not use a notation or any name for the totient function, but rather defined it in words. The first notation Euler used for the totient function appeared in one of his later works in 1784, where he used $\pi N$ to denote $\phi(N)$ and wrote ``$\pi1 = 0$" \cite{Euler_1775_pi_notation_for_phi}, yet this did not become the commonly used notation nowadays. In 1801, Carl Friedrich Gauss introduced the notation $\phi(n)$, along with the conviction that $\phi(1)=1$, in his number theory book \textit{Disquisitiones arithmeticae} \cite{Gauss_uses_Phi_notation}. Gauss accomplished his book by collecting works of mathematicians such as Euler, Fermat, Euclid, Diophantus, and many others, combining them with his own results in the field of number theory, which created this masterpiece that is still regarded as one of the greatest works of mathematical analysis. Due to the significant impacts of this book, the notation $\phi$ became the well know notation for the totient function. However, it was not until half a century later that $\phi$ was given the name ``Toteint Function". James Joseph Sylvester, founder of the American Journal of Mathematics, proposed the name ``totient" in his 1879 article \textit{On certain ternary cubic-form equations} \cite{name_totient_came_from}. Three different great mathematicians from different times, collectively established the definition, notation, and name for the totient function. In addition to its significance in number theory, $\phi(n)$ also showcases the beauty of cumulative contributions in mathematics and the timelessness of mathematical discovery across generations.\\

The earliest studies of applying the totient function multiple times, or the iteration of the totient function, began around 1929 by the works of Pillai \cite{Earliest_Studyof_Iteration_on_totient_Function}, which mainly focused on the number of times applying the totient function to get a number to 1. Erdös, Granvilie, Pomerance, and Spiro continued studies about this iteration of the totient function in 1990 \cite{termination_of_phi^k}. We discussed some of their results in the \textit{Related Works} section.

\section{Definitions and Properties}\label{gisec}
To understand the Phi function, we first come across the concept of relatively prime, or co-prime. Two numbers are relatively prime if they share no common factors greater than 1. Formally:
\begin{framed}
\begin{definition}[Relatively Prime]
\label{defRelativelyPrime}
For positive integers $a$ and $b$, we say $a$ is relatively prime to $b$ if
\[
gcd(a,b) = 1.
\]    
\end{definition}
\end{framed}
Here $gcd(a,b)$ denotes the greatest common divisor of $a$ and $b$, which is defined as the largest number in the set of common divisors shared by $a$ and $b$. Note that 1 is relatively prime to 1 by definition because $gcd(1, 1)= 1$. 
\begin{framed}
    \begin{definition}[Definition of $\phi$]
    \label{defPhi}  
    For $n\geq 1$, let $\phi(n)$ denote the number of positive integers less than or equal to $n$ that are relatively prime to $n$.
    \end{definition}
\end{framed}
For example, we present below a table of the first 10 numbers and its totient values.
\begin{table}[H]
    \caption{Values of $\phi(n)$ for The First 10 Positive Integers} 
    \label{tbl_10}
    \centering 
    \setlength{\tabcolsep}{20pt} 
    \begin{tabular}{c c c } 
    \hline\hline 
    $n$ & $\{ x \mid 1\leq x \leq n, gcd(x,n) = 1\}$ & Value of $\phi(n)$ \\ [0.5ex] 
    \hline 
    1 & \{ 1 \} & 1 \\ [1ex]
    2 & \{ 1 \} & 1 \\ [1ex]
    3 & \{ 1, 2 \} & 2 \\ [1ex]
    4 & \{ 1, 3 \} & 2 \\ [1ex]
    5 & \{ 1, 2, 3, 4 \} & 4 \\ [1ex]
    6 & \{ 1, 5 \} & 2 \\ [1ex]
    7 & \{ 1, 2, 3, 4, 5, 6 \} & 6 \\ [1ex]
    8 & \{ 1, 3, 5, 7 \} & 4 \\ [1ex]
    9 & \{ 1, 2, 4, 5, 7, 8 \} & 6 \\ [1ex]
    10 & \{ 1, 3, 7, 9 \} & 4 \\ [1ex]
    \hline 
    \end{tabular}
    \label{table:nonlin} 
    \end{table}
From the 10 values in \textit{Table \ref{tbl_10}}, several interesting questions arise: Why are most of these $\phi(n)$ even numbers? Why is $\phi(n)$ one less than $n$ for prime numbers? Why does $\phi(2)  \phi(5) = \phi(10)$, $\phi(2)\phi(3) = \phi(6)$, but $\phi(2) \phi(4) \neq \phi(8)$? The last question will be explained using the concept of multiplicative functions. According to the book \textit{An Introduction to Number Theory} by Harold Stark \cite{introduction_to_number_theory_book}, the following defines multiplicative functions:
\begin{framed}
    \begin{definition}[Multiplicative Functions]
    \label{defMultiplicative}  
    If the function $f(n)$ is defined for all positive integers $n$, then we say that $f(n)$ is multiplicative if for all pairs of relatively prime positive integers $a$ and $b$,
    \[
    f(ab) = f(a)  f(b)
    \]
    Furthermore, if this is true for all pairs of positive integers, no matter if they are relatively prime to each other or not, then we say $f(n)$ is completely multiplicative.
    \end{definition}
\end{framed}
This multiplicative property of functions can be quite useful, because we can break big numbers into its prime powers and deal with those powers separately. To approach the proof of Euler's totient function's multiplicativity, we will first establish the Linear Congruence Lemma and the Chinese Remainder Theorem.
\begin{framed}
     \begin{lemma}[Linear Congruences]
    \label{lemma_linear_cong} 
    For a linear congruence equation, where $a, b,$ and $m$ are constants, $x$ is the variable, consider the congruence:
    \[
        ax \equiv b \pmod{m}
    \]
    Let $gcd(a, m) =g$, if $g| b$, then there are exactly $g$ solutions to $x \pmod{m}$.
    \end{lemma}
\end{framed}
\begin{proof}
Since $g| b$, $ b = gc$ for some integer $c$. Because $g$ is the greatest common divisor of $m$ and $a$, by Bezout's Identity, there exists integers $r$ and $s$ such that 
\[g = ra + sm.\] Substitute this into $b=gc$ and applying the distributive property of multiplication gives
\[
b \equiv rac + smc \equiv rac \pmod{m}
\]
It is apparent that $x \equiv rc \pmod{m}$ would be a solution satisfying $ax \equiv b \pmod{m}$. With this first solution, a total of $g$ solutions could be generated. It could be shown that the set of numbers
\[
\{ x_{i} = x+ i \left( \frac{m}{g}\right) \mid 0\leq i \leq g-1\}
\]
are all solutions to $ax = b \pmod{m}$, by substituting in $x_i$ and showing that the left side of the equation equals $b\pmod{m}$.
\[
ax_i = ax + im\left(\frac{a}{g}\right),
\]
since $g$ is a factor of $a$, then $\frac{a}{g}$ must be an integer. $i$ is also and integer, so 
\[im\left( \frac{a}{g}\right) \equiv 0 \pmod{m}\]
Therefore:
\[
ax_i \equiv ax+0 \equiv b \pmod{m}
\]
The set of solutions $x_i$ has $g$ elements, because $0\leq i \leq g-1$. This completes the proof that there must be $g$ solutions to the linear Congruence equation.

\end{proof}
With the linear congruence established, we can prove the Chinese Remainder Theorem, a Theorem first discovered by the Chinese mathematician Sunzi in the 5th century \cite{sunzi_suanjing}.
\begin{framed}
     \begin{theorem}[Chinese Remainder Theorem]
    \label{thmChinese_Remainder}  
    For relatively Prime integers $m_1, m_2$, there exists one unique solution $x\pmod{m_1m_2}$ such that 
    \[
    \begin{cases}
    x \equiv a \pmod{m_1}\\
    x \equiv b \pmod{m_2}
    \end{cases}
    \]
    \end{theorem}
\end{framed}
\begin{proof}
    Suppose $x_1$ is a solution that satisfies the first equation $x_1 \equiv a_1 \pmod{m_1}$. Then
    \[
    x_1 = a_1 +m_1 y 
    \]
    for some integer $y$, which could be substituted into the second equation $x_1 \equiv a_2 \pmod{m_2}$, yielding
    \[
    a_1 + m_1y = a_2 \pmod{m_2}.
    \]
    Rearranging the above equation gives:
    \[
    m_1y = a_2 - a_1 \pmod{m_2}
    \]
    This is a linear congruence equation with respect to $y$. Since $gcd(m1, m2)=1$, and $1 \mid (a_2-a_1)$ because $a_1,a_2$ are both integers, so \textit{Lemma \ref{lemma_linear_cong}} guarantees a solution to $y$, thus a solution to $x_1$. Furthermore, this solution is a unique solution $\pmod{m_1m_2}$. If there exist another solution $x_2$ such that
    \[
    \begin{cases}
    x_2 \equiv a \pmod{m_1}\\
    x_2 \equiv b \pmod{m_2}
    \end{cases}
    \]
    then $x_1 -x_2 \equiv 0 \pmod{m1}$ and $x_1 -x_2 \equiv 0 \pmod{m2}$, which means 
    \[
    \begin{cases}
    m_1 \mid x_1 -x_2\\
    m_2 \mid x_1 -x_2
    \end{cases}
    \]
    Since both $m_1$ and $m_2$ are factors of $x_1 - x_2$, and $gcd(m_1,m_2)=1$, this indicates that 
    \[
    m_1m_2 \mid x_1 -x_2.
    \]
    In order to divide $m_1m_2$, $x_1$ and $x_2$ must have the same remainder when dividing $m_1m_2$, which is also evident in the derivation below:
    \[
    x_1 - x_2 \equiv 0 \pmod{m_1m_2}
    \]
    \begin{equation} \label{eq_1}
    x_1 \equiv x_2 \pmod{m_1m_2}
    \end{equation}
    
    \textit{Equation \ref{eq_1}} shows that $x_1$ and $x_2$ end up as the same solution$\pmod{m_1m_2}$, completing the proof that there exists only one unique solution$\pmod{m_1m_2}$ that satisfies the system of congruences.
\end{proof}

The Euler totient function $\phi(n)$ is a multiplicative function, but not completely multiplicative. First, it is easy to show that it is not completely multiplicative, because $2\cdot 4 =8$, but $\phi(2) \phi(4) =2 \neq \phi(8)$.
Second, we prove its multiplicativity:
\begin{framed}
    \begin{theorem}[Mutiplicativity of $\phi(n)$]
    \label{thmMultiplicative}  
    For positive integers $m$ and $n$, if $gcd(m,n) =1$, then
    \[
    \phi(mn) = \phi(m)  \phi(n).
    \]
    \end{theorem}
\end{framed}
\begin{proof}
Define three sets:
\[
    \Z^{'}_{m} =\{   k \mid 1\leq k \leq m,\space gcd(k,m) = 1\}\]
 \[
    \Z^{'}_{n} =\{   k \mid 1\leq k \leq n,\space gcd(k,n) = 1\}
    \]
    \[\Z^{'}_{mn} =\{   k \mid 1\leq k \leq mn,\space gcd(k,mn) = 1\}
\]
by \textit{Definition \ref{defPhi}}, the number of elements in $\Z^{'}_{n}$ is equal to $\phi(n)$, which is the same for other sets, so \[
    \left|\Z^{'}_{m} \right|=\phi(m)
    \]
\[
     \left|\Z^{'}_{n}\right| = \phi(n)
     \]\[
    \left|\Z^{'}_{mn} \right|=\phi(mn)
\]
Define the function $f(x) = (x\bmod{m} ,\space  x\bmod n)$,  $$f: \Z^{'}_{mn} \to \Z^{'}_{m} \times \Z^{'}_n.$$
For an order pair $(a,b)$ from the set $\Z^{'}_{m} \times \Z^{'}_{n}$, there is a corresponding $x$ in $\Z^{'}_{mn}$ such that 
\[
\begin{cases}
 x \equiv a \pmod{m}
\\ x \equiv b \pmod{m}
\end{cases}
\]
With the precondition of $gcd(m,n) = 1$, the \textit{Chinese Remainder Theorem} guarantees the existence of a solution $x$. For every output in the cross product set, there exists at least one input in the first set that produces the output, indicating that the function $f(x)$ is a surjection. So
$$\left| \Z^{'}_{mn} \right| \geq \left| \Z^{'}_{m} \times \Z^{'}_n\right|.$$
Furthermore, the \textit{Chinese Remainder Theorem} also guarantees that the solution $x$ is unique $\bmod mn$, which means that the function $f(x)$ is also an injection. So
$$\left| \Z^{'}_{mn} \right| \leq \left| \Z^{'}_{m} \times \Z^{'}_n\right|.$$
Since $f(x)$ is a bijection function, the number of elements must be the same for the set of inputs and outputs.
\begin{equation}\label{eq4}
\left| \Z^{'}_{mn} \right| = \left| \Z^{'}_{m} \times \Z^{'}_n\right|.
\end{equation}
By the \textit{Fundamental Principle of Counting}, the number of elements in $\Z^{'}_m \times \Z^{'}_n$ is equal to the product of the number of elements in $\Z^{'}_m$ and $\Z^{'}_n$, so $\left| \Z^{'}_{m} \times \Z^{'}_n\right| = \left| \Z^{'}_{m}\right| \left|\Z^{'}_n\right|$. When this is substituted into \textit{Equation \ref{eq4}},
\[
\left| \Z^{'}_{mn} \right| = \left| \Z^{'}_{m} \right|  \left| \Z^{'}_n\right|,
\]
which proves
\begin{equation}
\phi(mn) = \phi(m) \phi(n) \notag
\end{equation}
with the precondition of $gcd(m,n) =1$.

\end{proof}

To calculate $\phi(n)$ for large numbers of $n$, it is difficult to count all relatively prime numbers with $n$ like \textit{Table \ref{tbl_10}}. Fortunately, we can calculate $\phi(n)$ with the following formula, if its prime factors are known:
\begin{framed}
    \begin{theorem}[Formula to Calculate $\phi(n)$] By the Fundamental Theorem of Arithmetic, suppose the unique factorization of $n$ into powers of primes is that
$
n = p_1^{a_1} p_2^{a2}p_3^{a_3}\cdots p_k^{a_k},
$ where $p_i$ denotes distinct prime factors and $a_i \in \Z^+$. Then 
    \label{thmPhi}  
    \[ 
    \phi(n) = n \prod_{i=1}^{k} \left( 1- \frac{1}{p_i} \right).
    \]
    Note that the only exception to this formula is $\phi(1)$. By \textit{definition \ref{defRelativelyPrime}}, $\phi(1)=1$.
    \end{theorem}
\end{framed}
\begin{proof} The \textit{Fundamental Theorem of Arithmetic} guarantees the factorization of n, so
\[
\phi(n) = \phi(p_1^{a_1} p_2^{a_2}p_3^{a_3}\cdots p_k^{a_k}).
\]
Since the only prime factor of $p^c$ is $p$, $p_1^{a_1},p_2^{a_2},\cdots,p_k^{a_k}$ are pairwise relatively prime to each other, then the \textit{Multiplicative Property} of the totient function could be applied:
\[
\phi(n) = \prod_{i=1}^{k} \phi(p_i^{a_i}).
\]
For each $\phi(p_i^{a_i})$, by \textit{Definition \ref{defPhi}}, it is equal to number of positive integers in the set 
\[
S =\{ x \mid 1\leq x\leq\ p_i^{a_i}, \space gcd(x,p_i^{a_i})=1 \}.
\]

It is clear that a number of $p_i^{a_i}$ integers fulfill the first condition. Because the only prime factor of $p_i^{a_i}$ is $p_i$, the numbers that are not relatively prime to $p_i^{a_i}$ must be multiples of $p_i$. Likewise, only multiples of $p$ are relatively prime to $p_i^{a_i}$. Now, to count these numbers, a multiple of $p_i$ occurs every $p_i$th element in the set of integers that fulfills the first condition, so there are a total of $$\frac{p_i^{a_i}}{p_i} = p_i^{a_i -1}$$
elements that are not relatively prime to $p_i^{a_i}$. Removing these elements from the $p_i^{a_i}$ integers would obtain the set $S$. Then the number of elements in set $S$ is 
\[
\left| S \right| = p_i^{a_i} - p_i^{a_i-1} = p_i^{a_i}\left( 1 - \frac{1}{p_i}\right).
\]
Therefore, $\phi(p_i^{a_i}) = p_i^{a_i}\left( 1 - \frac{1}{p}\right)$, so 
\[
\phi(n) = \prod_{i=1}^{k} p_i^{a_i}\left( 1 - \frac{1}{p_i}\right).
\]
Bringing each $p_i^{a_i}$ out of the product and substituting in $n = p_1^{a_1} p_2^{a_2}p_3^{a_3}\cdots p_k^{a_k}$ gives:
\[
\phi(n) = n\prod_{i=1}^{k} \left( 1 - \frac{1}{p_i}\right),
\]
which completes the proof.
\end{proof}

With the formula, the following corollaries would follow:
\newenvironment{specialcorollary}[1]{
    \renewcommand{\thecorollary}{#1}\corollary
}{
    \endcorollary
}

\begin{framed}
    \begin{specialcorollary}{\ref{thmPhi}.1}[Constant Property] \label{corSame}
    If $a$ contains all prime factors of $b$, then
    \[
    \phi(a b) = \phi(a)  b
    \]
    \end{specialcorollary}
\end{framed}
\begin{proof}
Let $S_a$ denote the set of prime factors of $a$, $S_b$ denote the set of prime factors of $b$. Suppose $S_a = \{ p_a \mid 1\leq a \leq k\}$, with $p_a$ representing the prime factors of $a$. Because $a$ contains all prime factors of $b$, $$S_b \in S_a,$$
so the set of prime factors of $ab$ is still $S_a$. By \textbf{Theorem \ref{thmPhi}},
\[
\phi(ab) = ab \prod_{a=1}^{k}\left( 1-\frac{1}{p_a}\right)
\]
Again, by \textbf{Theorem \ref{thmPhi}},
$
\phi(a) = a \prod_{a=1}^{k}\left( 1-\frac{1}{p_a}\right)
$, then
\[
\phi(a)\cdot b = ab \prod_{a=1}^{k}\left( 1-\frac{1}{p_a}\right)
\]
As shown above, the left side and the right side of the Corollary statement equals the same expression, completing the proof.
\end{proof}
 \begin{framed}
    \begin{specialcorollary}{\ref{thmPhi}.2}[Strictly Decreasing Property]\label{cor_decrease}
    For all positive integers $n>1$, 
    \[
    \phi(n) < n
    \]
    \end{specialcorollary}
\end{framed}
\begin{proof}
Because $n\geq2$, $n$ must contain at least one prime factor. Let $\{ p_i \mid 1\leq i \leq k \}$ denote the distinct prime factors of $n$. $\frac{1}{p_i} > 0$ for all $p_i$,  so 
\[
0 <1-\frac{1}{p_i} < 1
\]
Multiplying $n$ by positive numbers smaller than one will always result in a smaller number, so
\begin{equation}
n\prod_{i=1}^{k} \left( 1 - \frac{1}{p_i}\right) <n \notag
\end{equation}
By \textit{Theorem \ref{thmPhi}}, the left side of the above inequality equals $\phi(n)$. Therefore,
\[
\phi(n)<n
\]
\end{proof}
These two Corollaries offer assistance in proving theorems about iterations of the totient function.

\section{Iterations of the Totient Function}
\begin{framed}
    \begin{theorem}[Convergence to 1 of Iterated Totient Function]
    \label{thmConvergeto_1}
    For any positive integer $n$, there exists an $m$ such that
   $$\underbrace{\phi(\cdots\phi(\phi(\phi(}_{\textit{$\phi()$ $m$ times}}n)))\cdots) = 1$$
    \end{theorem}
\end{framed}

\begin{proof}
    First of all, the range and the domain of the Euler's totient function are all positive integers, so the output of the function would always be a valid input to the function, which guarantees that no matter how many times the totient function is applied to a positive integer, it will always give a positive integer and it would never be undefined throughout this process. By \textit{Corollary \ref{cor_decrease}}, 
    $$\phi(n)<n,$$
    which indicates that the value of the expression decreases at least by 1 every time the totient function is applied. Then if the totient function is applied $n-1$ times to $n$, it must decrease to 1, the smallest possible value in the range, because even if it only decreases by 1 each time, $n - (n-1) =1.$
    \[
    \underbrace{\phi(\cdots\phi(\phi(\phi(}_{\textit{$\phi()$ $n-1$ times}}n)))\cdots) =1.
    \]
    If the totient function is applied more times, it would stay at 1 because $\phi(1) =1$. The iterations of the totient function converges to 1, so for any $m \geq n-1$,
    \[
    \underbrace{\phi(\cdots\phi(\phi(\phi(}_{\textit{$\phi()$ $m$ times}}n)))\cdots) = 1,
    \]
    which completes the proof.
\end{proof}

\textit{An example of the convergence to 1 is as follows:} \\
Let $n =5$, we record the iterations of the totient function applied to $n$ below,
\begin{align*}
    \phi(5) &=4 \\
    \phi(\phi(5)) &= \phi(4) = 2\\
    \phi(\phi(\phi(5))) & =\phi(2) =1 \\
    \phi(\phi(\phi(\phi(5)))) & =\phi(1) =1 \\
    \vdots \\
    \phi(\phi(\cdots \phi(5)\cdots)) & = \phi(1) = 1 
\end{align*}
It is evident that applying the totient function to the number 5 four times or more would equal 1. Taking this theorem a step further, not only would the iterated totient function lead any number to 1, but it also must equal 2 at some iteration before it converges to 1. To prove that, we first show the specialty of the numbers 2 and 1.

\begin{framed}
    \begin{lemma}[Odd and Even of $\phi(n)$]
    \label{lemmaEven}  
    Through observation, $\phi(2) = 1$ and $\phi(1)=1$, an odd number. For all other integers $n>2$, let $m = \phi(n)$, then 
    \[
        m \text{ must be even.}
    \]
    \end{lemma}
\end{framed}

\begin{proof}
    For simplicity, the problem can be broken into two cases:
    \begin{itemize}
    \item \textbf{Case 1} ($n$ only includes the prime factor 2)\textbf{.}\\
    Let $n = 2^a$ for some integer $a > 1$. Since 2 is the only prime factor of $n$, by \textit{Theorem \ref{thmPhi}}:
    \begin{equation}
        \phi(n) = n\cdot(1-\frac{1}{2}) = \frac{1}{2}\cdot n \notag
    \end{equation}
    Substituting in $n=2^a$, we have:
    \begin{equation}
        \phi(n) = \frac{1}{2}\cdot 2^a = 2^{a-1} = 2^{a-2}\cdot2 \notag
    \end{equation}
    Because $a \in \Z^+$, $a>1$ indicates that $a \geq 2$, which ensures $2^{a-2} \geq  1$. Now, $m = \phi(n)$ could be used,
    \begin{equation}
        m \equiv 2^{a-2}\cdot2 \equiv 2^{a-2}\cdot 0 \equiv 0 \pmod{2}. \notag
    \end{equation}
    \item \textbf{Case 2} ($n$ includes prime factors other than 2)\textbf{.}\\
    Let $n = p^{k}\cdot a$, with $p \nmid a$, and $p$ being any prime factor other than 2. $gcd(p^k,a)=1$, which means $p^k$ is relatively prime with a $a$ by \textit{Definition \ref{defRelativelyPrime}}, so the multiplicative property could be applied to $\phi(n)$ by \textit{Theorem \ref{thmMultiplicative}:}
    
    \begin{equation}
        \phi(n) = \phi(p^k\cdot a) = \phi(p^k)\cdot \phi(a).\notag
    \end{equation}
    
    To calculate $\phi(p^k)$, by \textbf{Theorem \ref{thmPhi}}:
    \begin{equation}
        \phi(p^k) = p^k \cdot (1 -\frac{1}{p}) = p^{k-1}\cdot (p - 1).\notag
    \end{equation}
    Now, substitute in $ m = \phi(n)$:
    \begin{equation}
        m = \phi(p^k)\cdot \phi(a) = \phi(a)\cdot p^{k-1}\cdot (p - 1). \notag
    \end{equation}
    Because $p$ is a prime number other than 2, $p-1$ must be an even number.    $k \geq 1$, which ensures that $p^{k-1} \geq 1$, along with $\phi(a)\geq 1$. Therefore:
    \begin{equation}
        m\equiv \phi(a)\cdot p^{k-1}\cdot 0 \equiv 0 \pmod{2} \notag
    \end{equation}
    \end{itemize}
In both cases, which together include all values of $n>2$, it is clear that $ m \equiv 0 \pmod{2}$. This completes the proof that $m$ must be an even number, except for the output 1 coming from special inputs 1 and 2.
\end{proof}

This Lemma leads to the theorem that for any positive integer $n$, when we apply the phi function to $n$ multiple times, such as $\phi(\cdots\phi(\phi(\phi(n)))\cdots)$, it would reach 2 at some point during the iteration.
\begin{framed}
    \begin{theorem}[Iteration of Totient Values to 2]
    \label{thmConverge}
    For any positive integer $n$ greater than 2, there exist some finite positive integer $m < n-1$ such that
   $$\underbrace{\phi(\cdots\phi(\phi(\phi(}_{\textit{$\phi()$ $m$ times}}n)))\cdots) = 2$$
    \end{theorem}
\end{framed}

\begin{proof}
From \textit{Theorem \ref{thmConvergeto_1}}, we know that iterating the totient function $n-1$ times or more to $n$ would give:
\begin{equation} \label{eq_10}
\underbrace{\phi(\cdots\phi(\phi(\phi(}_{\textit{$\phi()$ $n-1$ times}}n)))\cdots) =1.
\end{equation}
Let $(m+1)$ be the least number such that
\[
\underbrace{\phi(\cdots\phi(\phi(\phi(}_{\textit{$\phi()$ $(m+1)$ times}}n)))\cdots) = 1,
\]
By \textit{Equation \ref{eq_10}}, $(m+1) \leq n-1$, so $m < n-1$. Let $$x = \underbrace{\phi(\cdots\phi(\phi(\phi(}_{\textit{$\phi()$ $m$ times}}n)))\cdots).$$
Then 
\[
\phi(x) = 1.
\]
\textit{Lemma \ref{lemmaEven}} states that the only solutions to x are 
\[
x = 1 \quad\text{or} \quad x=2.
\]
However, $ x= 1 $ can be eliminated because $m+1$ is the first iteration such that the value of $\phi(\phi(\cdots \phi(n)\cdots))$ equals 1. $x \neq 1$, therefore, the only solution to $x$ is
\[
 x =2.
\]
This completes the proof that 
$$
\underbrace{\phi(\cdots\phi(\phi(\phi(}_{\textit{$\phi()$ $m$ times}}n)))\cdots) = 2
$$
for some finite integer $m < n-1$.
\end{proof}
\textit{An example of Theorem \ref{thmConverge} is as follows:} \\
Consider the number $27$. When we record values of the iterates of the totient function applied to 27,
\begin{align*}
    \phi(27) &= 18 \\
    \phi(\phi(27)) &= \phi(18) = 6\\
    \phi(\phi(\phi(27))) & =\phi(6) =2 \\
    \phi(\phi(\phi(\phi(27)))) & =\phi(2) =1 
\end{align*}
It is evident that at the third iteration, $\phi(\phi(\phi(27))) =2$. It could be tried for any arbitrary positive integer, and it will always reach the number 2 at some point when the totient function is iterated enough times. Furthermore, we notice an interesting property with the above example:
\[
1+2+6+18 = 27.
\]
When the sum of all the substeps until $\phi^{(4)}(27)=1$ is computed, it is exactly 27. This is not a coincidence. In fact, all numbers in the form of $3^k$ has this property.

\begin{framed}
    \begin{specialcorollary}{\ref{thm_Kth_Power}.1}[Sum of Iterated Totient Functions] \label{cor_3^k}
    For all $n=3^{k}$, where $k \in \Z^{+}$,
    \[
    n = \phi(n) + \phi(\phi(n))+ \phi(\phi(\phi(n))) + \cdots + \phi(2).
    \]
    The iteration is terminated at $\phi(2)$ because it is the first time $\phi(\phi(\cdots \phi(n)\cdots)) = 1$. Adding more 1s after that would not provide more insight. This result follows directly from \textit{Theorem \ref{thm_Kth_Power}}.
    \end{specialcorollary}
\end{framed}
\begin{proof}
We proceed by induction on the variable $k$.
  \begin{itemize}
    \item \textbf{Base Case} ($k=1$)\textbf{.}\\
    When $k=1$, we check that the statement holds true by simply calculating step by step. 
    \[
    \phi(3) = 2, \quad \phi(\phi(3)) =\phi(2) = 1.
    \]
    Then the statement holds true because 
    \[
    \phi(3) +\phi(\phi(3)) =2+1 = 3.
    \]
    \item \textbf{Induction Hypothesis} ($k=a$)\textbf{.}\\
    For some fixed but arbitrary integer $a \geq 1$, assume that the statement holds true for $k=a$. Therefore, 
    \[
    \phi(3^a) + \phi(\phi(3^a))+\phi(\phi(\phi(3^a))) + \cdots + \phi(2) = 3^a.
    \]
    \item \textbf{Induction Step} ($k=a+1$)\textbf{.}
    The goal is to show that
    \[
    \phi(3^{a+1}) + \phi(\phi(3^{a+1}))+\phi(\phi(\phi(3^{a+1}))) + \cdots + \phi(2) = 3^{a+1},
    \]
     which we establish by evaluating the left side of the above equation term by term. \textit{Theorem \ref{thmPhi}} tells that $$\phi(3^{a+1}) = 3^{a+1}\left( 1 - \frac{1}{3}\right) = 3^{a}\cdot 2.$$
     For the second term $\phi(\phi(3^{a+1}))$, $3^{a}$ is relatively prime to 2, so we can apply the \textit{multiplicative property}:
     $$\phi(\phi(3^{a+1})) = \phi(3^a)\cdot \phi(2) = \phi(3^a).$$
     The third term equals the totient value of the second term, so $\phi(\phi(\phi(3^{a+1}))) = \phi(\phi(3^a))$. By similar logic, all the following terms would be $$\underbrace{\phi(\cdots\phi(\phi(}_{\textit{$\phi()$ $m$ times}}3^{a+1}))\cdots) = \underbrace{\phi(\cdots\phi(\phi(}_{\textit{$\phi()$ $m-1$ times}}3^a))\cdots),$$ until the last term $\phi(2)$. Then 
     \begin{align}\notag
     &\phi(3^{a+1}) + \phi(\phi(3^{a+1}))+\phi(\phi(\phi(3^{a+1}))) + \cdots + \phi(2) \\
     \label{eq_16}
     & = 3^{a} \cdot 2 + \phi(3^a) + \phi(\phi(3^a)) + \cdots + \phi(2).\\
     \label{eq_17}
     & = 3^{a}\cdot 2  + 3^{a} \\
     & = 3^a (2+1) = 3^{a+1} \notag
     \end{align}
     \textit{Equation \ref{eq_16}} to \textit{Equation \ref{eq_17}} is a valid step that utilizes the \textit{Induction Hypothesis}. The statement is indeed true for $k = a+1$.
    \end{itemize}
We established a base case, and proved that if the statement is true when $k=a$ for any arbitrary number $a \geq 1$, it would imply the correctness of the statement for the next number $k = a+1$, which completes the proof by induction. 
\end{proof}

One may wonder why is the number 3 so special such that the sums of iterated totient values of its powers equals its original value. First, $3-1 = 2$, so $\phi(3-1)=1$, which ``goes away" neatly because 1 is trivial in multiplication. For other prime numbers, if $p-1 \neq 2$, then $\phi(p-1) \neq 1$. Therefore, when the phi function is applied again, both $ p-1 $ and $\phi(p-1)$ would be even by \textit{Lemma \ref{lemmaEven}}, so the multiplicative property could not be applied because the totient function is not a completely multiplicative function. However, there is a type of special number for which we found the general formula on the sums of iterated totient values -- Fermat Primes.
\begin{framed}
    \begin{definition}[Fermat Numbers]   
    Fermat numbers $F_{n}$ are numbers in the form of
    \[
    F_{n} = 2^{2^{n}}+1.
    \]
    \end{definition}
\end{framed}
\textit{Some examples of Fermat Numbers include:}
\begin{align*}
     3&= 2^{2^0} +1 , \quad \text{is prime.} \\
     65537 &= 2^{2^4} +1, \quad \text{is prime.} \\
     4294967297 & = 2^{2^5} +1, \quad \text{is composite.}
\end{align*}
The sequence of Fermat numbers can be found on \textit{The On-Line Encyclopedia of Integer Sequences}$:$ A000215 \cite{OEIS_sequence}. Fermat primes are the prime numbers in this sequence. It was originally conjectured by Fermat that all numbers in this sequence are prime. In 1732, Euler disproved Fermat's conjecture by showing that $F_5 = 4294967297$ is composite because it is divisible by $641$ \cite{Euler_disprove_Fermat_Primes}. The primality of all Fermat Numbers still remains a puzzle for modern mathematics.
\begin{framed}
    \begin{theorem}[Fermat Primes]  
    \label{thm_Fermat_Primes}
    If a prime number $p$ is in the form of $ p = 2^{k}+1$, where $k$ is a positive integer, then $k$ must only contain the prime factor $2$, which is:
    \[
     k = 2^{n}
    \]
    for some integer $n$.
    \end{theorem}
\end{framed}

\begin{proof}
    Suppose, for contradiction, that there exists a positive integer $k$ that  contains prime factors other than $2$, such that $2^{k} + 1$ is prime. Let $p_{1}$ denote the prime factor of $k$ that is not 2, then $k = p_{1}a$ for some positive integer $a$. Therefore:
    \begin{equation}
    2^{k} +1 = 2^{p_{1} a} +1 = (2^{a})^{p_{1}} + 1 \notag
    \end{equation}
    
    This equation is in the form of the following polynomial. Let
    \begin{equation}
    f(x) = x^{m} + 1, \text{where $m\in \Z^{+}$} \notag
    \end{equation}
    When $x=2^{a}$ and $ m=p_{1}$, $f(x) = 2^k+1$. When $m$ is an odd number, $x = -1$ would be a root to the polynomial $f(x)$ because $(-1)^{m} + 1 =0$. Therefore, $(x+1)$ would be a factor to the polynomial $f(x)$, meaning that there exist some polynomial $g(x)$ such that 
    \begin{equation}
    f(x) = (x+1)\cdot g(x) \notag
    \end{equation}
    Because $p_1$ is a prime factor other than 2, $p_1$ is an odd number since the only even prime is 2. Then, if $m=p_1$, then $(x+1)$ would be a factor of $f(x)$ as shown above. Substituting in $x=2^{a}$, $(2^a + 1)$ would be a factor of $2^k+1$:
    \begin{equation}
    (2^a+1) \vert  2^k+1 \notag
    \end{equation}
    Since $a\in \Z^{+}$, it is easy to show that $2^a+1 \geq 3$. Since $p_1 > 2$, $2^{a}+1 <(2^{a})^{p_1}+1 = 2^k+1 $. This indicates that $2^k+1$ has a factor $2^{a}+1$ that is not 1 and not itself, contradicting the supposition that $2^k+1$ is a prime number.
    This proves our assumption false. Thus, there does not exist a positive integer $k$ that contains prime factors other than $2$, such that $2^{k} + 1$ is prime.
\end{proof}
Since all prime numbers in the form of $2^k+1$ must be a Fermat prime, we would use $2^{k}+1$ to denote Fermat primes for simplicity.

\begin{framed}
    \begin{lemma}[Sum of Iterated Totient on Fermat Primes]
    \label{lemmaFermat}  
    For all Fermat primes, primes of the form $p = 2^{k} +1$ where $k \in \Z^{+}$,
    \[
    \phi(p) + \phi(\phi(p))+ \phi(\phi(\phi(p)))+ \cdots + \phi(2) = 2p -3 
    \]
    \end{lemma}
\end{framed}

\begin{proof}
    Because $p$ is a prime number, by \textit{Theorem \ref{thmPhi}}, the first term $\phi(p) = p-1$. Then the second term would be $$\phi(p-1) = \phi(2^k) = 2^{k-1}.$$
    Then the third term would be $\phi(\phi(\phi(p)))=\phi(2^{k-1}) = 2^{k-2}$. Because $2$ is the only prime factor involved, all following terms would form a \textit{Geometric Sequence} that decrease by a common ratio of 2 until the last term equals 1, which is $\phi(2)$ in the original equation.
    \begin{align*}
        \phi(p) + \phi(\phi(p))+ \phi(\phi(\phi(p)))+ \cdots + \phi(2) = p-1 + 2^{k-1} +2^{k-2} + \cdots + 2 +1.
    \end{align*}
    To evaluate the sum of the above geometric sequence, we apply the ``Rolling the Snowball" method. First, add 1 and subtract 1 from the equation, then $1+1 =2$, which could be added to 2, so $2+2 =4$, and it keeps going on as a snowball rolling bigger and bigger, eventually reaching the last term $2^{k-1}+2^{k-1} = 2^k$.
    \begin{align*}
        2^{k-1} +2^{k-2} + \cdots + 2 +1 +1 -1 &= 2^{k-1}+2^{k-2}+\cdots +4+2+2 -1 \\
        &=2^{k-1}+2^{k-2}+\cdots +8+4+4 -1 \\
        & = 2^{k-1}+2^{k-2}+\cdots+8+8 -1 \\
        & \vdots \\
        & = 2^{k-1}+2^{k-1} -1\\
        & = 2^{k} -1
    \end{align*}
    $2^{k} -1 = p -2$, so putting everything together gives:
    \[
    \phi(p) + \phi(\phi(p))+ \phi(\phi(\phi(p)))+ \cdots + \phi(2) = p-1 + p-2 = 2p -3,
    \]
    which completes the proof of the sum of iterated totient functions on Fermat Primes.
\end{proof}
This result complies with the sum of iterated totient values of $3^k$, because $2\cdot 3 -3=3$. Similar to powers of three, iterated totient sums on Fermat Primes could also be generalized to the $k$th power. Using results from \textit{Lemma \ref{lemmaFermat}}, we found the formula for $p^2$ is:
\[
 \phi(p^2) + \phi(\phi(p^2))+ \phi(\phi(\phi(p^2)))+ \cdots + \phi(2) =2p^2 -3p.
\]
The pattern seems obvious, with 2 and 3 being constant coefficients, and the power of $p$ increasing. Therefore, we soon conjectured that:
\[
 \phi(p^k) + \phi(\phi(p^k))+ \phi(\phi(\phi(p^k)))+ \cdots + \phi(2) =2p^k -3p^{k-1}.
\]
We checked that the above equation agrees with the proven formula for $3^k$, because $2\cdot 3^k - 3\cdot 3^{k-1} = 3^k$. Everything points in favor of this formula, but unfortunately, it is not the correct formula. When we observe a pattern, such as Fermat's Primes conjecture and the above formula, no matter how obvious and appealing the pattern is, we could never establish it as the truth before proving it. Modern computers emerge fast, but despite the power to check billions of numbers for a pattern, it is still not regarded as true for all numbers. Rigorous proofs are required to establish formulas, using the same methods ancient Greeks proved their math.\\

The general formula of Fermat Primes to the $k$th power is as follows:
\begin{framed}
    \begin{theorem}[Fermat Primes to $k$th Power]
    \label{thm_Kth_Power}
    Let $p$ be a Fermat prime, For all $p^k$ where $k \in \Z^{+}$,
    \[
    \phi(p^k) + \phi(\phi(p^k))+ \phi(\phi(\phi(p^k))) + \cdots + \phi(2) = \frac{2}{p+1}\left[p^k(p-1)+2\left(\frac{p-1}{2}\right)^k\right] - 1
    \]
    \end{theorem}
\end{framed}
\begin{proof}
    We utilize induction on the variable $k$.
    \begin{itemize}
    \item \textbf{Base Case} ($k=1$)\textbf{.}\\
    As proved above in \textit{Lemma \ref{lemmaFermat}}, the left side of the equation at $k=1$ equals $2p-3$.
    \begin{equation}
     \phi(p) + \phi(\phi(p))+ \cdots + \phi(2) = 2p -3. \notag
    \end{equation}
    Substituting $k=1$ to the right side of the statement gives
    \begin{equation}
     \frac{2}{p+1}\left[p(p-1)+2\left(\frac{p-1}{2}\right)\right]-1= \frac{2}{p+1}(p^2-p+p-1) -1. \notag
    \end{equation}
    By the difference of squares, $p^2-1 = (p-1)(p+1)$, so
    \begin{equation}
     \frac{2}{p+1}(p^2-1) -1 = 2(p-1) -1 = 2p-3. \notag
    \end{equation}
    It is shown that both the left side and right side of the equation equals $2p-3$, which verifies the truth of the statement for $k=1$.
    \item \textbf{Induction Hypothesis} ($k=a$)\textbf{.}\\
    For some fixed but arbitrary integer $a \geq 1$, assume that the statement holds true for $k=a$:
     \begin{equation}
     \phi(p^a) + \phi(\phi(p^a))+ \cdots + \phi(2) = \frac{2}{p+1}\left[p^a(p-1)+2\left(\frac{p-1}{2}\right)^a\right] - 1. \notag
    \end{equation}
    \item \textbf{Induction Step} ($k=a+1$)\textbf{.}\\ 
    The goal is to prove:
    \begin{equation}
     \phi(p^{a+1}) + \phi(\phi(p^{a+1}))+ \cdots + \phi(2) = \frac{2}{p+1}\left[p^{a+1}(p-1)+2\left(\frac{p-1}{2}\right)^{a+1}\right] - 1. \notag
    \end{equation}
    To achieve this, we analyze the left side of the equation term by term. The only prime factor of $p^{a+1}$ is $p$, so by \textit{Theorem \ref{thmPhi}}, 
    \begin{equation}
     \phi(p^{a+1})  = p^{a+1}\cdot\left(1-\frac{1}{p}\right) = p^a\cdot(p-1) . \notag
    \end{equation}
    For the second term, $\phi(\phi(p^{a+1})) = \phi(p^a \cdot (p-1))$. The only prime factor of $p^a$ is $p$, and clearly $p$ is not a factor of $p-1$, so $gcd(p^a,p-1) = 1$. By \textit{Definition \ref{defRelativelyPrime}}, $p^a$ is relatively prime with $p-1$. Then the multiplicative property of the totient function holds, with \textit{Theorem \ref{thmMultiplicative}} yielding:
    \begin{equation}
     \phi(\phi(p^{a+1})) = \phi[p^a \cdot (p-1)] = \phi(p^a)\cdot \phi(p-1). \notag
    \end{equation}
    For the third term, $\phi(\phi(\phi(p^{a+1})))= \phi[\phi(p^a)\cdot \phi(p-1)]$. It is known that $\phi(p^a) = p^{a-1}\cdot (p-1)$. Because $p$ is a Fermat Prime: $p-1 = 2^n$, and $\phi(p-1) = 2^{n-1}$. The only prime factor present in $\phi(p-1)$ is $2$, which is also present in $\phi(p^a)$, then by \textit{Corollary \ref{corSame}},
    \begin{equation}
     \phi(\phi(\phi(p^{a+1})))= \phi[\phi(p^a)\cdot \phi(p-1)] = \phi[\phi(p^a)]\cdot \phi(p-1). \notag
    \end{equation}
    For the fourth term, by similar logic, $\phi(\phi(p^a))$ also contains the factor 2, so we can then apply \textit{Corollary \ref{corSame}} again,
    \begin{equation}
     \phi(\phi(\phi(\phi(p^{a+1}))))= \phi\left[ \phi(\phi(p^a))\cdot \phi(p-1)\right] = \phi\left[\phi(\phi(p^a))\right] \cdot \phi(p-1).\notag
    \end{equation}

    For the following terms, $\phi(\phi(\cdots \phi(p^a)\cdots ))$ will always contain the factor 2, so we can continue applying \textit{Corollary \ref{corSame}} to all following terms, $$ \phi(\underbrace{\phi[\phi(\cdots \phi(}_{\textit{i}}p^a)\cdots )]) = \underbrace{\phi[\phi(\cdots \phi(}_{\textit{i-1}}p^a)\cdots )] \cdot \phi(p-1) ,$$until the term $m+1$, where 
    \begin{equation}\underbrace{\phi(\cdots\phi(\phi(\phi(}_{\textit{$\phi()$ $m$ times}}p^a)))\cdots) = \phi(2). \notag
    \end{equation} 
     The term $m+1$ is $\phi(2)\cdot \phi(p-1)$, which is just $\phi(p-1)$. The $(m+2)$th term would then be $\phi(\phi(p-1))$, and all following terms would be $\phi(\phi(\cdots \phi(p-1)\cdots))$. Again by \textit{Theorem \ref{thmConverge}}, this would eventually converge to 2, leading to our last two terms:

    \begin{equation}
        \phi(\phi(\cdots \phi(p-1)\cdots)) = 2, \notag
    \end{equation}
        \text{and the last term } 
    \begin{equation}
        \phi(2)= 1. \notag
    \end{equation}
    For clarity, we recorded the values of the above derivations in the table below:

    \begin{table}[H]
    \caption{Left Side of Equation} 
    \centering 
    \setlength{\tabcolsep}{20pt} 
    \begin{tabular}{c c c } 
    \hline\hline 
    Term\# & Expression in Equation & Equivalent Value \\ [0.5ex] 
    \hline 
    1 & $\phi(p^{a+1})$ & $p^a \cdot (p-1)$\\ [2ex]
    2 & $\phi(\phi(p^{a+1}))$ & $\phi(p^a) \cdot \phi(p-1)$ \\ [2ex]
    3 & $\phi(\phi(\phi(p^{a+1})))$ & $\phi(\phi(p^a)) \cdot \phi(p-1)$ \\ [2ex]
    $\vdots$ & $\vdots$ & $\vdots$ \\ [2ex]
    $i$ & $\underbrace{\phi(\phi(\cdots\phi}_{\textit{$\phi()$ $i$ times}}(p^{a+1})\cdots))$ & $\underbrace{\phi(\phi(\cdots\phi}_{\textit{$\phi()$ $i-1$  times}}(p^{a})\cdots)) \cdot \phi(p-1)$ \\ [2ex]
    $\vdots$ & $\vdots$ & $\vdots$  \\ [2ex] 
    $m+1$ & $\underbrace{\phi(\phi(\cdots\phi}_{\textit{$\phi()$ $m+1$ times}}(p^{a+1})\cdots))$ & $\phi(2)\cdot \phi(p-1)$ \\ [4ex]
    $m+2$ & $\phi(\phi(\cdots\phi(p^{a+1})\cdots))$ & $\phi(\phi(p-1))$ \\ [2ex]
    $m+3$ & $\phi(\phi(\cdots \phi(p^{a+1})\cdots))$ & $\phi(\phi(\phi(p-1)))$ \\ [2ex]
    $\vdots$ & $\vdots$ & $\vdots$ \\ [2ex]
    $m+i$ & $\underbrace{\phi(\phi(\cdots\phi}_{\textit{$\phi()$ $m+i$ times}}(p^{a+1})\cdots))$ & $\underbrace{\phi(\phi(\cdots\phi}_{\textit{$\phi()$ $i$ times}}(p-1)\cdots))$ \\ [2ex]
    $\vdots$ & $\vdots$ & $\vdots$ \\ [2ex]
    second to last & $\phi(\phi(\cdots \phi(p^{a+1})\cdots))$ & $2$ \\ [2ex]
    last & $\phi(\phi(\cdots \phi(p^{a+1})\cdots))$ & $\phi(2)$ \\ [2ex]
    \hline 
    \end{tabular}
    \label{tbl_2} 
    \end{table}
    
    To Add up all terms in \textit{Table \ref{tbl_2}}, $\phi(p-1)$ can be factored out for terms from 2 to $m+1$.
    \begin{align} \notag
    \phi(&p^{a+1}) + \phi(\phi(p^{a+1}))+ \phi(\phi(\phi(p^{a+1})))+\cdots + \phi(2)  \\ \notag
    &= p^a\cdot (p-1) + \phi(p-1)\cdot [\phi(p^a)+\phi(\phi(p^a)) + \cdots +\phi(2)]\\ \label{eq34}
    &\quad +\phi(\phi(p-1)) + \phi(\phi(\phi(p-1))) + \cdots + \phi(2).
    \end{align}
    Recall that \textit{Lemma \ref{lemmaFermat}} states:
    \begin{equation}
        \phi(p)+\phi(\phi(p))+\phi(\phi(\phi(p)))+\cdots + \phi(2) = 2p-3.  \notag
    \end{equation}
    Substituting $ \phi(p-1) =\phi(\phi(p)) $ into the above equation gives:
    \begin{equation}
        \phi(p)+\phi(p-1)+\phi(\phi(p-1))+\cdots + \phi(2) = 2p-3.  \notag
    \end{equation}
Rearranging this equations gives:
    \begin{align} \notag
        \phi(\phi(p-1))+\phi(\phi(\phi(p-1))) + \cdots + \phi(2) &= 2p-3 -\phi(p) -\phi(p-1)\\ 
        \label{equA}
        & = p-2 - \phi(p-1)\\ 
        \label{equB}
        & =\phi(p-1) -1
    \end{align}
    \textit{Equation \ref{equA}} to \textit{Equation \ref{equB}} is valid because $p-1 = 2 \cdot \phi(p-1)$.
    Now, using the above for the sum of terms $m+2$ to the last term, along with the \textit{Induction Hypothesis} for the sum of terms 2 to $m+1$, \textit{Equation \ref{eq34}} becomes:
    $$\phi(p^{a+1}) + \phi(\phi(p^{a+1}))+ \cdots + \phi(2)
    = p^a\cdot (p-1) + \phi(p-1)\cdot (\frac{2}{p+1}[p^a(p-1)+2(\frac{p-1}{2})^a] - 1)$$
    \begin{equation}
    + \phi(p-1) -1 \notag
    \end{equation}
    Notice that $-1$ in the parenthesis multiplies $\phi(p-1)$, which nicely cancels out with the $+\phi(p-1)$ on the second row, a reason for leaving the $\phi(p-1)$ notation. Then, with some algebraic manipulations:
    \begin{align*}
        \phi(p^{a+1}) + \phi(\phi(p^{a+1}))+ \cdots + \phi(2) &= p^a\cdot(p-1)+ \frac{p-1}{p+1} [ p^a(p-1)+ 2(\frac{p-1}{2})^a] -1\\
        &= \frac{p-1}{p+1} [p^a(p+1) + p^a(p-1)+2(\frac{p-1}{2})^a]-1\\ 
        &= \frac{p-1}{p+1} [2\cdot p^{a+1}+2(\frac{p-1}{2})^a]-1\\
        &= \frac{2}{p+1}[p^{a+1}(p-1)+2(\frac{p-1}{2})^{a+1}] -1 
    \end{align*}
    This shows that the statement holds true for $k=a+1$.

    \end{itemize}
We demonstrated that if the statement holds for $k=a$, it would also hold for the next number $k=a+1$, so with the base case, it follows that for any positive integer $k$, \textit{Theorem \ref{thm_Kth_Power}} applies.
\end{proof}

The uniqueness of Fermat Primes lies in the fact that $2$ is the only even prime number. When the totient function is applied to primes in the form of $p^k +1$, it becomes $p^k$, which makes the recursive sum easy to track because $p^k$ only contains one prime factor. However, the primes in the form $p^k +1$ must be a Fermat Prime, because a prime number other than 2 is odd, making $p^k+1$ an even number thus divisible by 2. If $p=2$, \textit{Theorem \ref{thm_Fermat_Primes}} states that it must be a Fermat Prime.

\section{Related Works}
\subsection{Termination of $\phi^k(n)$}
In the works of Erdös, Granvilie, Pomerance, and Spiro, \textit{On the Normal Behavior of the Iterates Of some Arithmetic Functions}, they proved some theorems regarding iterations of the totient function, and stated some open problems \cite{termination_of_phi^k}. Their study was heavily focused on applying the totient function $k$ times to a number $n$, such that it is the first time $\phi^k(n) = 1$, which means the iteration is terminated at the first iteration when the iterated totient value equals 1. They defined the function $k(n) = k$ as the least number such that $\phi^k(n) = 1$. Some key results include: If $n = 2^j$, then,
\[
k(n) = j = \frac{\log n}{\log2}.
\]
We could apply this to \textit{Lemma \ref{lemmaFermat}}, so that not only we know the sum would be $2p-3$, but also we could know that there will be a total of 
$$k(p) =1+\frac{\log (p-1)}{\log 2}$$ 
terms in the sum of totient values series. Another result they found is if $n = 2\cdot 3^j$, then
\[
k(n) = j+1 = \left\lceil \frac{\log n}{\log 3} \right\rceil,
\]
where $\lceil x \rceil$ denotes the least integer greater or equal to $x$. We could apply this to \textit{Corollary \ref{cor_3^k}} to figure out the total number of terms in that series. Furthermore, the case $n=2^j$ and $n = 2\cdot 3^j$ is actually the two extreme behaviors of $k(n)$ \cite{termination_of_phi^k}, so 
\[
\left\lceil \frac{\log n}{\log 3} \right\rceil \leq k(n) \leq \left\lceil \frac{\log n}{\log2} \right\rceil,
\]
for any positive integer $n$. Further investigations could determine the total number of terms in \textit{Theorem \ref{thm_Kth_Power}}, the series of totient values on Fermat Primes to the $k$th power. \\

\textit{An Arithmetic Function Arising From the $\phi$ Function}, an article from Princeton University by Harold Shapiro, also studied the termination of the iterated totient function \cite{C(n)_function}. However, his definition of the last iterate is a little different. He established similar results as \textit{Theorem \ref{thmConverge}} --- the iterated totient value always arrives at 2 --- so he defined the function $C(n) = x$ such that $\phi^x(n) = 2$. One key result they produced with strong relation to this paper is that
\[
C(2^x +1) = x  \quad \text{ if and only if $2^x+1$ is a Fermat Prime}.
\]
Combined together with \textit{Theorem \ref{thm_Kth_Power}}, these results might assist in the ultimate investigations concerning the primality and compositeness of all Fermat Numbers.

\subsection{Sum of Totients}
There are also other approaches to summing up totient functions, and interesting theorems related. For example, in the book \textit{History of the Theory of Numbers} by Leonard Eugene Dickson \cite{History_of_the_theory_of_numbers_6/pi^2}, they present a theorem by summing up totient values from 1 to $n$:
\[
\lim_{n \to \infty} \frac{2}{n^2} \sum_{k=1}^{n} \phi(k) = \frac{6}{\pi^2}.
\]
Here an interesting connection arise. Proved by Euler in his \textit{Desummis Serierum Reciprocarum}, the infinite series $\frac{1}{n^2} $ converges to $\frac{\pi^2}{6}$ \cite{series_1/n^2_pi^2/6},
\[
 \lim_{n \to \infty} \sum_{k=1}^{n} \frac{1}{k^2} = \frac{\pi^2}{6},
\]
which is the reciprocal of the totient sums. Although these two topics are far from related, the results show a strong connection between the totient sums and the $\frac{1}{n^2}$ series. Future investigations could try to uncover the reason behind such interesting connection

\section{Conclusion}
This paper proposes closed-form formulas for the sums of iterated totient functions. First, A broad picture of number theory and functions were introduced, leading to the Euler’s totient Function. We explained the name and purpose of the function, introducing iteration approaches on the totient function. Next, we presented definitions and fundamental theorems including relatively prime, multiplicative functions, the Linear congruence Lemma, and the Chinese remainder theorem. We prove the formula to calculate $\phi(n)$, with two corollaries following. We used the second corollary immediately in the proof of totient function’s convergence, while we used the first corollary multiple times in the proof of the iterated totient sums of Fermat Primes, the main result of this paper:
$$\phi(p^k) + \phi(\phi(p^k))+ \phi(\phi(\phi(p^k)))+ \cdots + \phi(2) = \frac{2}{p+1}\left[p^k(p-1)+2\left(\frac{p-1}{2}\right)^k\right] - 1.$$
This theorem produces a clean corollary for the first Fermat Prime, the number 3:
$$ \phi(3^k) + \phi(\phi(3^k))+\phi(\phi(\phi(3^k))) + \cdots + \phi(2) = 3^k.$$
Other related works in the field of the iterated totient function provides a method to track the number of terms in these series, as the sums of iterated totient functions can be expressed in the series form:
$$\sum_{i=1}^{k\left(3^n\right)} \phi^i(3^n) = 3^n.$$
Overall, the approach of iteration sums opens up a new topic to investigate. Many theorems and generalizations await to be uncovered.
\newpage
\bibliography{refrences}

@misc{introduction_to_number_theory_book,
  title={An Introduction to Number Theory},
  author={{Harold M. Stark}},
  year={1987},
  howpublished={The MIT Press},
}

@misc{euler_totient,
  title={Theoremata arithmetica nova methodo demonstrata},
  author={{L. Euler}},
  year={1763},
  howpublished={Commentarii academiae scientiarum Petropolitanae, Volume 8, pp. 141-146},
}

@misc{Euler_disprove_Fermat_Primes,
  title={Observationes de theoremate quodam Fermatiano aliisque ad numeros primos spectantibus},
  author={{L. Euler}},
  year={1738},
  howpublished={Commentarii academiae scientiarum Petropolitanae, Volume 6, pp. 103-107},
}

@misc{series_1/n^2_pi^2/6,
  title={De Summis Serierum Reciprocarum},
  author={{L. Euler}},
  year={1740},
  howpublished={Commentarii academiae scientiarum Petropolitanae, Volume 7, pp. 123-134},
}

@misc{Euler_1775_pi_notation_for_phi,
  title={Speculationes circa quasdam insignes proprietates numerorum},
  author={{L. Euler}},
  year={1784},
  howpublished={Acta Academiae Scientiarum Imperialis Petropolitanae, Volume 1780: II, pp. 18-30.},
}

@misc{chinese_remainder_thm,
  title={Chinese Remainder Theorem: Applications In Computing, Coding, Cryptography},
  author={{Dingyi Pei, Arto Salomaa, Cunsheng Ding}},
  year={1996 Oct. 25},
  howpublished={World Scientific}}

@misc{sunzi_suanjing,
  title={Sunzi Suanjing (The Mathematical Classic of Master Sun)},
  author={{Sunzi}},
  year={2007},
  howpublished={Translated by  J.N.Crossley and A.W.C.Lun},
}

@misc{OEIS_sequence,
   title={The OEIS Foundation},
  author={{N.J.A.Sloane}},
  year={2016},
  howpublished={The On-Line Encyclopedia of Integer Sequences, published electronically at \url{http://oeis.org}},
}

@misc{termination_of_phi^k,
  title={On the Normal Behavior of the Iterates Of some Arithmetic Functions},
  author={{Erdös, P., Granvilie, A., Pomerance, C., Spiro, C. }},
  year={1990},
  howpublished={ In: Berndt, B.C., Diamond, H.G., Halberstam, H., Hildebrand, A. (eds) Analytic Number Theory. Progress in Mathematics, vol 85. Birkhäuser Boston. \url{https://doi.org/10.1007/978-1-4612-3464-7_13}},
}

@misc{History_of_the_theory_of_numbers_6/pi^2,
  title={History of the Theory of numbers},
  author={{Leonard Eugene Dickson }},
  year={1919},
  howpublished={Published by the Carnegie Institution of Washington},
}

@misc{Earliest_Studyof_Iteration_on_Totient_Function,
  title={On a function connected with $\phi(n)$},
  author={{S. S. Pillai}},
  year={1929},
  howpublished={Bull. Amer. Math. Soc. 35 (6) 837 - 841},
}

@misc{Gauss_uses_Phi_notation,
  title={Disquisitiones Arithmeticae},
  author={{Carl Friedrich Gauss}},
  year={1801},
  howpublished={Königliche Gesellschaft der Wissenschaften zu Göttingen Vol. 1.},
}

@misc{Totient_in_Latin,
  title={The Words of Mathematics},
  author={{Schwartzman, Steven}},
  year={1994},
  howpublished={Washington, DC: Mathematical Association of America},
}

@misc{name_totient_came_from,
  title={On Certain Ternary Cubic-Form Equations},
  author={{J. J. Sylvester}},
  year={1879},
  howpublished={American Journal of Mathematics, Vol. 2, No. 4, pp. 357-393},
}

@misc{Open_Conjecture_Carmichael,
  title={On Euler’s $\phi-$function},
  author={{Carmichael, Robert Daniel}},
  year={1907},
  howpublished={Bulletin of the American Mathematical Society 13, 241-243 },
}

@misc{Open_Conjecture_Lehmer,
  title={On Euler's totient function},
  author={{D. H. Lehmer }},
  year={1932},
  howpublished={Bulletin of the American Mathematical Society Vol.38},
}

@misc{RSA,
  title={A method for obtaining digital signatures and public-key cryptosystems},
  author={{Rivest, Ronald L., Adi Shamir, and Leonard Adleman}},
  year={1978},
  howpublished={Communications of the ACM 21.2, 120-126.},
}

@misc{Number_theory_oldest_branches_in_mathematics,
  title={Number Theory: An Introduction to Pure and Applied Mathematics},
  author={{Redmond, D}},
  year={1978},
  howpublished={ CRC Press. \url{https://doi.org/10.1201/9781003067535}},
}
\bibliographystyle{plain}
\end{document}